\documentclass[draft,12pt,oneside,intlimits]{article}
\usepackage{amsmath}
\usepackage{amsthm}
\usepackage{amssymb}
\usepackage{enumerate}

\def\CC{\mathbb{C}}
\def\QQ{\mathbb{Q}}

\def\NN{\mathbb{N}}
\def\RR{\mathbb{R}}

\def\one{\mathbf{1}}

\DeclareMathOperator*{\dd}{{d\hskip-1pt}}

\newcommand*{\bR}{{\mathbb R }}
\def\RRb{\overline{\mathbb R}}

\newcommand*{\bN}{{\mathbb N }}

\newcommand*{\fM}{{\mathfrak{M}} }
\newcommand*{\fMe}{\ensuremath{ {\mathfrak{M}}_1  }}

\def\sstar{\#}
\def\conv{\mathop{\mathrm{conv}}}
\def\cconv{\mathop{\mathrm{\overline{conv}}}}

\newcommand*{\ws}{weak$^*$}

\theoremstyle{plain}
\newtheorem{proposition}{Proposition}
\newtheorem{theorem}[proposition]{Theorem}
\newtheorem{lemma}[proposition]{Lemma}
\newtheorem{corollary}[proposition]{Corollary}

\theoremstyle{definition}
\newtheorem{definition}[proposition]{Definition}

\theoremstyle{definition}
\newtheorem{example}[proposition]{Example}
\newtheorem{remark}[proposition]{Remark}

\DeclareMathOperator*{\supp}{supp}

\begin{document}

\hyphenpenalty=10000 \hfuzz=10pt \brokenpenalty=10000
 \numberwithin{proposition}{section}
\numberwithin{equation}{section}

\title{Potential theoretic approach to rendezvous numbers}

\author{B\'{a}lint Farkas, Szil\'{a}rd Gy.~R\'{e}v\'{e}sz}
\date{}

\maketitle
\begin{center}
\small Alfr\'{e}d R\'{e}nyi Institute of Mathematics\\
Hungarian Academy of Sciences\\ Re\'{a}ltanoda utca 13--15, H-1053, Budapest, Hungary\\
revesz@renyi.hu\\[3ex]
Technische Universit\"{a}t Darmstadt\\ Fachbereich Mathematik, AG4 \\ Schlo{\ss}gartenstra\ss
e 7, D-64289 Darmstadt, Germany\\
farkas@mathematik.tu-darmstadt.de

\end{center}
\vskip2ex
\begin{abstract}
We analyze relations between various forms of energies (reciprocal
capacities), the transfinite diameter, various Chebyshev constants
and the so-called rendezvous or average number. The latter is
originally defined for compact connected metric spaces $(X,d)$ as
the (in this case unique) nonnegative real number $r$ with the
property that for arbitrary finite point systems
$\{x_1,\dots,x_n\}\subset X$, there exists some point $x\in X$ with
the average of the distances $d(x,x_j)$ being exactly $r$. Existence
of such a miraculous number has fascinated many people; its
normalized version was even named ``the magic number'' of the metric
space. Exploring related notions of general potential theory, as set
up, e.g., in the fundamental works of Fuglede and Ohtsuka, we
present an alternative, potential theoretic approach to rendezvous
numbers.

 \vskip2em \noindent{\small
\emph{AMS Subj.~Clas.(2000):} Primary: 31C15. Secondary: 28A12, 54D45} \vskip1em
\end{abstract}

\let\oldthefootnote\thefootnote\def\thefootnote{}
\footnotetext{Supported by OTKA (Hungarian Scientific Research
Fund)}
\let\thefootnote\oldthefootnote
\section{Introduction}\label{sec:Intro}

Rendezvous numbers, studied usually in metric spaces, have already
attracted much attention. The existence and uniqueness of such a
miraculous number in compact\newline connected metric spaces was
shown by O.~Gross \cite{Gross}. Later many authors were fascinated
by the topic. They calculated the rendezvous numbers in particular
cases, and extended the results in the direction of weak rendezvous
numbers or rendezvous numbers of unit spheres in Banach spaces (see,
e.g., \cite{baronti/casini/papini:2000}, \cite{Bj}, \cite{CMY},
\cite{GVV}, \cite{Lin}, \cite{MN}, \cite{NY}, \cite{Stadje},
\cite{Th}, \cite{W1,W2,W3} and \cite{YZ}).

Already Bj\"{o}rck applied certain tools of potential theory in studying
constants related to rendezvous numbers \cite{Bj}. Now, exploring
notions of general, abstract potential theory, as set up, e.g., in
the fundamental works of Fuglede and Ohtsuka -- and, in particular,
drawing from some exploration of mutual energies, Chebyshev
constants and transfinite diameters over locally compact topological
spaces with lower semicontinuous, nonnegative and symmetric kernels
(see \cite{Fug2}, \cite{Oh2}, \cite{Bela}) -- we arrive at an
understanding of these quantities from a more general viewpoint.

However, to achieve this, we need to recover and even partially extend the relevant basic
material. In particular, we thoroughly investigate energies and Chebyshev constants, and even
also their ``minimax duals'' in function of \emph{two} sets. The technical reason for that is
that the classical definitions are kind of saddle point special cases of these more general
notions, and we need to utilize special monotonicity and other properties, which stay hidden
when considering only the diagonal cases.

Let us recall the appropriate setting of potential theory in locally compact spaces. First,
$+\infty$ is added to the set of real numbers, i.e., we let
$\overline{\bR}:=\bR\cup\{+\infty\}$ endowed with its natural topology such that
$\overline{\bR}_+$ will be compact. Moreover, we will use the notation $\conv E$ for the convex
hull of a subset $E\subset \RR$ and $\cconv E$ for the closed convex hull in
$\overline{\RR}_+$, meaning, for example, $\overline{\conv}(0,+\infty)=[0,+\infty]$.

Throughout the paper $X$ denotes a locally compact Hausdorff
space, and $k:X\times X\rightarrow \RRb$ is a kernel function in
the sense of Fuglede \cite[p.~149]{Fug}. That is, we assume that
$k$ is {\em lower semicontinuous} (l.s.c.) as a two variable
function over $X\times X$, and that $-\infty < k(x,y) \le
+\infty$. Moreover, in this paper we assume that $k\ge 0$, and
that $k$ is symmetric, i.e., $k(x,y)=k(y,x)$ for all $x,y\in X$.

Under these standard assumptions there exist \emph{energy-minimizing} (or, as sometimes called,
\emph{capacitary}) measures (see, e.g., \cite[Theorem 2.3]{Fug}), and also the equivalence of
the energies (reciprocal capacities) follows, see \cite[p.~159]{Fug}.

The vector space of all finite, signed regular Borel measures on
$X$ is $\fM^{\pm}(X):=\fM(X)-\fM(X)$, where
$$
\mathfrak{M}:=\mathfrak{M}(X):=\{\mu : \mu \, \mbox{is a
positive, regular Borel measure on} \, X, \, ||\mu||<\infty \}.
$$

We will also denote by $\mu$ the Carath\`{e}odory extension of $\mu$ to an outer measure, while
$\mu$-measurability of a set $H$ refers to the usual Carath\`{e}odory measurability.  Also, let
$\fMe(X)$ be the set of probability measures from $\fM(X)$, $\fMe:=\fMe(X):= \{\mu \in \fM(X) :
\mu(X) =1\}$. We say that $\mu$ is \emph{concentrated} on a set $H\subset X$, if each compact
set intersects the complement $X\setminus H$  of $H$ in a set of zero (outer-) measure. In our
case all measures are finite, hence this is the same as requiring measurability and having full
measure, cf.~\cite[p.~146]{Fug}. Furthermore, for an arbitrary set $H\subset X$, let
$\fMe(H):=\{\mu \in \fMe :\:\mbox{$\mu$ is concentrated on $H$}\}$. We also say that $\mu$ is
\emph{supported} in $H$ if $\supp \mu$ is a subset of $H$. Note that, if $H$ is a closed set,
then the two notions -- concentrated and supported on $H$ -- are the same. For convenience we
introduce the notation $\fMe^{\sstar}(H)$ for finitely supported probability measures on $H$.

The customary topology on $\fM$ (or on $\fM^{\pm}$) is the vague
topology which is the locally convex topology determined by the
seminorms $\mu\mapsto \left|\int_X f \dd\mu\right|$, $f\in
C_c(X)$, where $C_c(X)$ denotes the space of continuous functions
with compact support. In most places we consider only the family
$\fMe(K)$ of probability measures supported on the same compact
set $K$. In this case, by the Riesz Representation Theorem,
$\fM^{\pm}(K)=C(K)'$, moreover, the weak$^*$-topology determined
by $C(K)$ and the vague topology coincide. We will use nets to
describe convergence in $X$ or in $\fM$. For example, we will use
the filtering family $\mathcal{K}$ (directed upwards) of all
compact subsets of $X$. Given $\mu\in\fM(X)$ we define its traces
$\mu_K$ on compact sets\footnote{Here and also throughout the
sequel, the notation $K\Subset H$ means that $K$ is a compact
subset of $H$.} $K\Subset X$ by $\mu_K(A):=\mu(A\cap K)$. Using
the above terminology the regularity of $\mu$ implies that
$\lim_{\mathcal{K}} \mu_K=\mu$ in the vague topology. This, among
others, means that $\lim_{\mathcal{K}}\|\mu_K\|=\|\mu\|$.

\vskip1em

The potential and energy of $\mu\in \fM$ are defined as
\[
    U^\mu(x) := \int_X  k(x,y)\dd\mu(y) ~, \qquad
    W({\mu}) := \int_X\int_X  k(x,y)\dd\mu(y)\dd\mu(x)  \ .
\]
The existence of these two integrals follow because $k\ge 0$ and
it is lower semicontinuous -- however, $W({\mu})$ and $U^\mu$ may
attain $+\infty$ as well.

There are three different definitions of energy (reciprocal
capacity) in use in the fundamental work of Fuglede, see
\cite[p.~150]{Fug}, \cite[Eq.~(1) and (2), p.~153]{Fug}. All these
are related to the supremum of the potential $U^\mu(x)$ over
particular subsets of $X$.

In some less widely known works, however the theory have already been extended to generalized
mutual energies and ``dual energies'' of \emph{two sets} as variables, see \cite{Fug2} and
\cite{Oh2}: moreover, even generalized, two-variate Chebyshev constants were defined in
\cite{Oh2} as well. We add here the notion of dual Chebyshev constants, too. These allow to
present a fairly general approach to rendezvous numbers.

We try to make this paper as self-contained, as possible. In particular, we use only a basic
familiarity with the monographic paper of Fuglede \cite{Fug}, which we consider and use as our
basic reference to potential theory in general. For all material above those directly quoted
from \cite{Fug} we give a precise description, or at least sketch a proof. Nevertheless, for
historical completeness we will refer to various potential theoretical papers and many
potential theoretical results throughout.

\section{Definitions and preliminaries}
For our purposes we consider the following notion of ``energy'', among whose special cases will
be the three quantities used in \cite[p.~150--153]{Fug} (and recalled also in Definition
\ref{def:Fugledeenergy} below).

\begin{definition}\label{def:qq}
Let $H\subset X$ be fixed, and $\mu\in\fMe(X)$ be arbitrary. First
put
\begin{equation}\label{d:QmuH}
Q(\mu;H):=\sup_{x\in H} U^{\mu}(x)~,\quad\text{and also}\qquad
\underline{Q}(\mu;H):=\inf_{x\in H} U^{\mu}(x)~.
\end{equation}
The \emph{quasi-uniform energy} and the \emph{restricted
quasi-uniform energy} of $H$ are
\begin{equation}\label{qcapacity}
q(H):= \inf_{\mu\in\fMe(H)} Q(\mu;H) \quad \text{and} \qquad
q^{\sstar}(H):= \inf_{\mbox{\scriptsize
$\mu\in\fMe^{\sstar}(H)$}} Q(\mu;H) ~.
\end{equation}
Moreover, we define analogously for any {\it two} sets
$H,L\subset X$ the quantities
\begin{alignat}{3}\label{qstarcapacity}
q(H,L)&:= \inf_{\mu\in\fMe(H)} Q(\mu;L) &&\quad \text{and} \qquad& q^{\sstar}(H,L)&:=
\inf_{\mbox{\scriptsize $\mu\in\fMe^{\sstar}(H)$}} Q(\mu;L) ~. \intertext{Furthermore, let us
recall the dual notions \cite{Oh2}} \label{qstarul}
\underline{q}(H,L)&:=\sup_{\nu\in\fMe(H)}\underline{Q}(\nu;L)&&\quad\mbox{and}\qquad&\underline{q}^{\sstar}(H,L)&:=\sup_{\mbox{\scriptsize$\nu\in\fMe^{\sstar}(H)$}}\underline{Q}(\nu;L)~,
\intertext{and, similarly to \eqref{qcapacity}} \label{qulcapacity} \underline{q}(H)&:=
\underline{q}(H,H) &&\quad \text{and} \qquad&\underline{q}^{\sstar}(H)&:=
\underline{q}^{\sstar}(H,H) ~.
\end{alignat}
\end{definition}
These extended definitions -- with different notation -- appear already in the work of Fuglede
\cite{Fug2} and Ohtsuka \cite{Oh2}.  Recall the following univariate versions of energies,
commonly used in linear potential theory.

\begin{definition}\label{def:Fugledeenergy}
\begin{enumerate}[a)]
\item Let $\mu\in\fMe$ be any measure. Then we write
\begin{align}\label{mucapacities}
U(\mu):&=\sup_{x\in X} U^{\mu}(x) = \sup_{K\Subset X} U(\mu_K)~,  \\
V(\mu):&=\sup_{x\in \supp \mu} U^{\mu}(x) = \sup_{K\Subset X} V(\mu_K)~,
\end{align}
the second forms, as well as an analogous statement for $W(\mu)$, following from \cite[Lemma
2.2.2 and Footnote 1 on p.~152]{Fug}.

\item Accordingly, the ``uniform'', ``de la Vall\'{e}e-Poussin'' and ``Wiener'' energies
(reciprocal capacities) of any set $H\subset X$ are
\begin{align}
u(H):=&\inf_{\mu\in\fMe(H)} U(\mu)=\inf_{\mu\in\fMe(H),~ \supp
    \mu \Subset H} U(\mu)=\inf_{K\Subset H} u(K)~, \label{ucapacity} \\
v(H):=&\inf_{\mu\in\fMe(H)} V(\mu)=\inf_{\mu\in\fMe(H),~ \supp
    \mu \Subset H} V(\mu)= \inf_{K\Subset H} v(K)~,  \label{vcapacity} \\
w(H):=&\inf_{\mu\in\fMe(H)} W(\mu)= \inf_{\mu\in\fMe(H),~ \supp \mu \Subset H} W(\mu) =
    \inf_{K\Subset H} w(K)~ \label{wcapacity},
\end{align}
respectively, where equivalence of the last forms can be proved based on \cite{Fug}. Instead of
doing that, we will prove the similar, but more general statement for $q$ below in Lemma
\ref{l:q}.
\end{enumerate}
\end{definition}

\begin{lemma}\label{lem:lsc}Let $H,L\subset X$. Equipping the set $\fM_1(H)$ with the vague topology
the following functions are lower semicontinuous
\begin{enumerate}[a)]
\item \hfil$\displaystyle\fM_1(H)\times L\ni(\mu,x)\mapsto
U^{\mu}(x)\left(=\int_X k(x,y)\dd\mu(y)\right)$\hfill
\item \hfil$\displaystyle\fM_1(H)\times\fMe(L)\ni(\mu,\nu)\mapsto W(\mu,\nu):=\int_X\int_X k(x,y)
\dd\mu(y)\dd\nu(x)$\hfill
\item \hfil$\displaystyle\fMe(H)\ni\mu\mapsto Q(\mu;L)\left( =\sup_{x\in L} U^{\mu}(x)\right)$\hfill
\end{enumerate}
\end{lemma}
\begin{proof}
Part a) and b) are taken from \cite[Lemma 2.2.1]{Fug}, while part
c) follows from part a) by noticing that $Q(\mu;L)$ is a
pointwise supremum of l.s.c.~functions.
\end{proof}

\begin{lemma}\label{l:q} For any two sets $H,L \subset X$ we have
\begin{align}
q(H,L)&= \inf_{\mbox{\scriptsize $\mu\in\fMe(H)$} \atop
\mbox{\scriptsize $\supp \mu \Subset H$}}Q(\mu;L)=\inf_{K\Subset H} q(K,L)  \label{qlemma}~, \\
q^{\sstar}(H,L) &= \inf_{\mbox{\scriptsize
$\mu\in\fMe^{\sstar}(H)$}} Q(\mu;L)=\inf_{K\Subset H}
q^{\sstar}(K,L) \label{qstarlemma}~.
\end{align}
\end{lemma}
\begin{proof}
First,
\begin{equation}
q(H,L):=\inf_{\mbox{\scriptsize $\mu\in\fMe(H)$}} Q(\mu;L) \le
\inf_{\mbox{\scriptsize $\mu\in\fMe(H)$}\atop\mbox{\scriptsize
$\supp \mu \Subset H$} } Q(\mu;L) = \inf_{K \Subset H} q(K;L)~
\label{qstarproof}
\end{equation}
is obvious. Second, for any measure $\mu\in\fMe(H)$ and compact
set $K\Subset H$ we have $Q(\mu;L)\ge Q(\mu_K;L)=||\mu_K||
Q(\nu_{\mu,K};L)$ with $\nu_{\mu,K}:=\mu_K/||\mu_K||$ now
satisfying $\nu_{\mu,K} \in \fMe(K)$. Plainly, for any measure
$\mu\in\fMe(X)$ we have $\mu_K\le \mu$, and thus for nonnegative
kernels $k$ also the potentials satisfy $U^{\mu_K}\leq U^{\mu}$.
Now let us select some sequence of increasing compact sets
$K_n\Subset H$ with $\mu(K_n)>1-\frac{1}{n}$. Such a sequence
exists by regularity of the measure $\mu$, since $\mu\in \fMe(H)$
entails that the set $H$ is $\mu$-measurable
(cf.~\cite[p.~146]{Fug}). Thus we find
\begin{equation}\label{qstarconverse}
\inf_{K \Subset H} q(K,L)\leq q(K_n,L)\leq Q(\nu_{\mu,K_n};L)\leq\frac{1}{1-1/n}Q
(\mu_{K_n};L)\leq\frac{1}{1-1/n}Q (\mu;L),
\end{equation}
for all $n\geq 1$ and hence $\inf_{K \Subset H} q(K,L)\leq Q
(\mu;L)$. Since here $\mu\in \fMe(H)$ is arbitrary, we obtain
$\inf_{K \Subset H} q(K,L)\leq q (H,L)$. On combining this with
\eqref{qstarproof} gives the last two formulations of
\eqref{qlemma}. The proof of \eqref{qstarlemma} is similar.
\end{proof}
\begin{remark}\label{whytwosets} Although one would like to have
$q(H)=\inf_{K\Subset H} q(K)$, this is false in general.
\end{remark}
\begin{example}
As an example one can consider a two point metric space $X=\{a,b\}$ with the discrete metric as
the kernel $k$, and take $H=X$. It is obvious that for $\# K =1$ we have $q(K)=0$, while
$Q(\mu;H)=\max (\mu(a),\mu(b))\ge 1/2$ for any $\mu\in \fMe$, hence $q(H)=1/2\gneqq
\inf_{K\Subset H} q(K) =0$. The reason for the occurrence of this difficulty is lack of
monotonicity of $q(H)$. On the other hand, fixing one variable of the two-set function $q$, the
functions $q(\cdot,L)$ and $q(H,\cdot)$ are monotonous (see also Lemma \ref{monotoneset}
below). That is why involving two-set functions is necessary here and throughout the paper.
\end{example}
Not surprisingly, a kind of a dual statement  holds for
$\underline{q}$, at least for compact $L$.
\begin{lemma}\label{l:-q-} Let $H\subset $X be arbitrary and $L\Subset X$ be compact. Then we have
\begin{align}
\underline{q}(H,L)&= \sup_{\mbox{\scriptsize $\mu\in\fMe(H)$}
\atop
\mbox{\scriptsize $\supp \mu \Subset H$}}\underline{Q}(\mu;L)=\sup_{K\Subset H} \underline{q}(K,L)  \label{-q-lemma}~, \\
\underline{q}^{\sstar}(H,L) &= \sup_{\mbox{\scriptsize
$\mu\in\fMe^{\sstar}(H)$}\atop \mbox{\scriptsize $\supp \mu
\Subset H$}}\underline{Q}(\mu;L)=\sup_{K\Subset H}
\underline{q}^{\sstar}(K,L) \label{-q-starlemma}~.
\end{align}
\end{lemma}
\begin{proof} The inequality
\begin{equation}
\underline{q}(H,L):=\sup_{\mbox{\scriptsize $\mu\in\fMe(H)$}}
\underline{Q}(\mu;L) \ge \sup_{\mbox{\scriptsize
$\mu\in\fMe(H)$}\atop\mbox{\scriptsize $\supp \mu \Subset H$} }
\underline{Q}(\mu;L) = \sup_{K \Subset H} \underline{q}(K;L)~
\label{-q-starproof}
\end{equation}
is trivial. Similarly to the proof of Lemma \ref{l:q} let
$\mu\in\fMe(H)$ be arbitrary and consider its traces $\mu_K$ on
compact sets $K$. Further define
$\nu_K:=\nu_{\mu,K}:=\mu_K/||\mu_K||$. (Note that $\|\mu_K\| > 0$
can be assumed by fixing some compact set $\widetilde{K}\Subset H$
with $\mu(\widetilde{K}) >0$ and considering only compact sets
$K\supset \widetilde{K}$ here.) Then clearly $\nu_K\in\fMe(K)$. By
lower semicontinuity of $U^{\nu_K}$ (cf.~Lemma \ref{lem:lsc} a))
and referring to compactness of $L$ we find an element $y_K\in L$
minimizing $U^{\nu_K}(y)$ on $L$. Because $\mu$ is a regular Borel
measure $\lim_{K\in\mathcal{K}}\nu_{K}=\mu$ in the vague topology.

Also, by compactness of $L$ we have a subnet $\mathcal N$ of
$\mathcal{K}$ such that $\lim_{K\in \mathcal{N}}y_K=y_{0}$ for
some $y_0$. Again by l.s.c.~we obtain (using $\nu_{K} \to \mu$
and $y_{K}\to y_0$ along $\mathcal N$, and referring to Lemma
\ref{lem:lsc} a)),
$$
\liminf_{\mathcal N} U^{\nu_{K}}(y_{K}) \geq U^\mu(y_0)\geq
\inf_{y\in L} U^\mu(y)=\underline{Q}(\mu;L).
$$
Summing up
$$
\sup_{K\Subset H}\underline{q}(K,L)\geq \liminf_{\mathcal N}
\inf_L U^{\nu_{K}} =\liminf_{\mathcal N} U^{\nu_{K}}(y_{K}) \geq
\underline{Q}(\mu;L)~.
$$
Taking supremum in $\mu$ and comparing to \eqref{-q-starproof}
yields the assertion. The proof of \eqref{-q-starlemma} is
analogous.
\end{proof}
The quantities of $n^{\rm th}$ diameters, transfinite diameters and Chebyshev constants were
generalized from the classical logarithmic kernel case to some more general kernels already by
P\'{o}lya, Szeg\H o \cite{Sz}, Carleson \cite{C} and Choquet \cite{ChoquetSem} the way we present
below. (Note that another direction of generalization, due to Zaharjuta \cite{Z}, is also
considered in $\CC^n$.)

\begin{definition}
Let $H\subset X$ be fixed. The \emph{$n^{\text{th}}$ diameter} of
$H$ is defined as follows.
\begin{align}\label{Ddef}
   D_n(H) := & \inf_{w_1,\ldots,w_n \in H}
    {\frac{2}{(n-1)\,n}}\Bigg(
    \sum_{1\le j < l \le n} k(w_j,w_l)
   \Bigg)~.
\end{align}
\end{definition}

\begin{definition} For an arbitrary $H \subset X$ the \emph{$n^{\text{th}}$ Chebyshev constant} of
$H$ is defined as
\[
     M_n(H) := \sup_{w_1,\ldots,w_n \in H} \inf_{x\in H}
    {\frac{1}{n}}\Bigg(
   \sum_{k=1}^{n} k(x,w_k)
   \Bigg)
   \ .
\]
\end{definition}
In fact, even in the classical literature another variant of the
Chebyshev constant occurs. Namely, for an arbitrary $H \subset X$
the \emph{modified $n^{\text{th}}$ Chebyshev constant} of $H$ is
defined as
\[
 C_n(H) := \sup_{w_1,\ldots,w_n \in X} \inf_{x\in H}
    {\frac{1}{n}}\Bigg( \sum_{k=1}^{n} k(x,w_k) \Bigg) \ ,
\]
that is, allowing the ``zeroes'' $w_j$ spread out in $X$, but
considering the values only on $H$.

To put these Chebyshev constants into a general framework was probably a reason why Ohtsuka
considered the following notion \cite{Oh2}. Our motivation is mentioned already in Remark
\ref{whytwosets} and will be even clearer in view of Proposition \ref{prop:dchc}.

\begin{definition} For arbitrary $H,L \subset X$ the
\emph{(general) $n^{\text{th}}$ Chebyshev constant} of $L$ with
respect to $H$ is defined as
\[
M_n(H,L):= \sup_{w_1,\ldots,w_n \in H} \inf_{x\in L}
    {\frac{1}{n}}\Bigg( \sum_{k=1}^{n} k(x,w_k) \Bigg).
\]
\end{definition}
 Remark that Ohtsuka \cite{Oh2} defines these quantities essentially without any
assumption on the kernel $k$. To prove the convergence of these sequences one has to assume
that the kernel does not take the value $-\infty$, then on compact set an l.s.c.~kernel will be
bounded from below. The symmetry assumption -- here -- does not play any role.

Finally, let us define -- in a slightly more general setting, that
is, forgetting about the metric, usually involved in the context
-- the (weak) rendezvous number(s), or average distance number(s)
of the space $X$, or even of subsets of $X$. Again, for good
reasons we define these notions in dependence of two sets as
variables.

\begin{definition}\label{def:randi}For arbitrary subsets $H, L \subset X$ the $n^{\text{th}}$
\emph{(weak) rendezvous set} of $L$ with respect to $H$ is
\begin{equation}\label{nthrandi}
R_n(H,L):=\bigcap_{w_1,\dots,w_n\in H} \cconv
\Bigl\{p_n(x):=\frac 1n \sum_{j=1}^n k(x,w_j)~~:~~ x\in L
\Bigr\}~.
\end{equation}
Correspondingly, one defines
\begin{alignat}{3}
\label{randidef} R(H,L)&:=\bigcap_{n=1}^\infty R_n(H,L)~,&\quad&&
R(H)&:=R(H,H)~. \intertext{Similarly, one defines the
\emph{(weak) average set} of $L$ with respect to $H$ as}
\label{avridef} A(H,L)&:=\bigcap_{\mu\in\fMe(H)} \cconv
\Bigl\{U^{\mu}(x)~~:~~ x\in L \Bigr\} ~,&\quad&& A(H)&:=A(H,H)~.
\end{alignat}
\end{definition}

\begin{remark}\label{AmH} Denoting the interval
\begin{equation}\label{Amuint}
A(\mu,L):=[\underline{Q}(\mu;L), Q(\mu;L)]=
\cconv\{U^{\mu}(x)~:~x\in L\}~,
\end{equation}
we see that $R_n(H,L)$, $R(H,L)$ and $A(H,L)$ are all of the form
$\bigcap_{\mu} A(\mu,H)$, with $\mu$ ranging over all averages of
$n$ Dirac measures at points of $H$, over all probability measures
finitely supported in $H$ (and of rational probabilities, but
compare to Lemma \ref{lem:raci} below) and over all of $\fMe(H)$,
respectively.
\end{remark}

\begin{remark}\label{sdc} If $k$ is a continuous kernel -- in particular
when it is a metric on $X$, -- then it suffices to take convex
hull instead of closed convex hull whenever $L$ is compact, since
then together with $k$ also $U^{\mu}(x)$ is continuous for any
probability measure $\mu$. Thus for compact subsets $L$ of metric
spaces a real number $r\in\RR_{+}$ belongs to $R(H,L)$ if and
only if for any finite system of (not necessarily distinct)
points $x_1,\dots,x_n \in H$ (of number $n\in\NN$ taken
arbitrarily) we always have points $y,z\in L$ satisfying
\begin{equation}\label{nthweakrendezvous}
\frac 1n \sum_{j=1}^n k(y,x_j)\le r \qquad\text{and}\qquad \frac
1n \sum_{j=1}^n k(z,x_j) \ge r~,
\end{equation}
which is the usual definition of weak rendezvous numbers in
metric spaces (see \cite{Th}). Moreover, in case the set $L$ is
connected, this is further equivalent to the existence of a
``rendezvous point'' $x\in L$ with
\begin{equation}\label{nthstrongrendezvous}
\frac 1n \sum_{j=1}^n k(x,x_j) = r~.
\end{equation}
In particular, for compact and connected $L$ in a metric space
(or in a locally compact space with continuous kernel $k$) the
rendezvous set $R(H,L)$ consists of a \emph{unique} point, say
$R(H,L)=\{r(H,L)\}$, if this latter property is satisfied only
for $r=r(H,L)$.
\end{remark}

\begin{remark}\label{sd} If $k$ is only l.s.c., also potentials are l.s.c.,
which entails that they take their infimum over compact sets.
Thus for compact $L$ the first half of the above equivalent
formulation \eqref{nthweakrendezvous} remains valid even for
general kernels. However, for the second part we must already
write that ``$\forall s<r ~ \exists z\in L$ such that $\frac 1n
\sum_{j=1}^n k(z,x_j) > s$''. Such modification of the
formulation is necessary also when we consider sets $L\subset X$
which are not compact, or when we are discussing the case when
$+\infty\in R(H,L)$. Clearly, in our settings $R_n(H,L)$,
$R(H,L)$ and $A(H,L)$ are subsets of $[0,\infty]$, but note that
traditionally rendezvous numbers or average numbers are
considered only among the reals. Hence even in metric spaces our
notions slightly differ from the usual ones regarding the role of
$+\infty$.
\end{remark}

\begin{example}
For example we say that even $+\infty\in R(H,L)$ if for all finite systems of (not necessarily
distinct) points $x_1,\dots, x_n \in H$ and for all real $s$, however large, there is $z\in L$
satisfying $\frac 1n \sum_{j=1}^n k(z,x_j) > s$. Thus, e.g., taking $X:=\RR$ with the usual
Euclidean metric $k(x,y):=|x-y|$, the four possible variations with the sets $\RR$ and
$I:=[0,1]$ yield $R(\RR,I)=\emptyset$, $R(I,\RR)=[1/2,\infty]$, $R(I):=R(I,I)=\{1/2\}$ (and
thus $r(I)=1/2$ exists uniquely; see \cite{CMY}) and $R(\RR):=R(\RR,\RR)=\{+\infty\}$, see
Proposition \ref{prop:reni}. Note that this last case classically would be interpreted as a
case when there is no rendezvous number, while in our notation this is a case of uniqueness
with $r(\RR)=+\infty$. Similar phenomena occur also in case of finite accumulation points not
belonging to the ``restricted sense rendezvous sets'', interpreted without closure.
\end{example}

The aim of the present work is to study properties of these sets
and set functions. First we look at the various energies
(reciprocal capacities) and Chebyshev constants associated to a
set $H \subset X$.

\section{Basic properties of the quantities defined}\label{sec:properties}

\label{sec:Cheby} Most of the results in this section were already obtained by Fuglede
\cite{Fug2} and Ohtsuka \cite{Oh2}. To be self-contained we present full proofs to them.

Let us start by recalling from \cite{ChoquetSem} and \cite{Oh2} the definition of the Chebyshev
constant and the transfinite diameter as the limit of the respective sequences. It is shown in
\cite{ChoquetSem} that the sequence of $n^{\text{th}}$ diameters is monotonically increasing
(cf.~\cite{Fek} or \cite{PSz} for the classical case; see also \cite{C, Bela}). The limit is
denoted by $D(H):=\lim_{n\rightarrow\infty}D_n(H)$ and is called the {\em transfinite diameter}
of $H$.

To prove the convergence of the sequence of $n^{\text{th}}$ Chebyshev constants the following
lemma can be used. (See \cite[p.\ 233]{Fek}. See also \cite{PSz}, part I, Ch.\ 3, \S1, Problem
98, p.\ 23 \& p.\ 198 or \cite{Bela}.)
\begin{lemma}[\bf Fekete\rm]\label{l:quasi} Let $(s_n)$ be a {\em quasi-monotonous} sequence of real
numbers, meaning that either $(n+m)s_{n+m}\le n s_n + m s_m$ (in
this case we say that $(s_n)$ is quasi-monotonically decreasing)
or $(n+m)s_{n+m}\ge n s_n + m s_m$ (the sequence $(s_n)$ is
quasi-monotonically increasing). Then $\lim_{n\to\infty} s_n =
\inf s_n$ (when $(s_n)$ is quasi-monotonically decreasing) or
$\lim_{n\to\infty} s_n = \sup s_n$ (when $(s_n)$ is
quasi-monotonically increasing), the $\inf$ ($\sup$) being either
finite or $-\infty$ ($+\infty$).
\end{lemma}

\begin{proposition}
For any $H,L\subset X$, the Chebyshev constants $M_n(H,L)$
converge, more precisely
\begin{equation*}
\sup_{n\in\NN} M_n(H,L)=\lim_{n\rightarrow \infty} M_n(H,L).
\end{equation*}
In particular, $\sup\limits_{n\in\NN}
M_n(H,H)=\lim\limits_{n\rightarrow \infty} M_n(H,H)$ and
$\sup\limits_{n\in\NN} M_n(X,H)=\lim\limits_{n\rightarrow \infty}
M_n(X,H)$.
\end{proposition}
\begin{proof}
Since $M_n(H,L)$ is quasi-monotone increasing, once $M_n(H,L)=+\infty$ holds, the subsequent
terms are all infinity as well, i.e., $M_{n+1}(H,L)=+\infty$. Hence $M_n(H,L)$ converges (in
the extended sense) to $+\infty$. If $M_n(H,L)$ is finite for all $n\in\NN$ Lemma \ref{l:quasi}
applies. For the classical Chebyshev constants themselves see \cite[p.\ 233]{Fek} or
\cite{PSz}, part I, Ch.\ 3, \S1, Problem 98, p.\ 23 \& p.\ 198 or \cite{Bela}. For more general
kernels and of a similar notation (with respect of taking logarithms of classical formulations)
see \cite[Chapter IV]{C}.
\end{proof}

The limit $M(H,L):=\lim_{n\rightarrow\infty}M_n(H,L)$ is called
the \emph{general Chebyshev constant of $L$ relative to $H$}. In
the two special cases the terminology for the limit is the
following: $M(H):=\lim_{n\rightarrow\infty}M_n(H)$ is called the
\emph{Chebyshev constant} of $H$ and
$C(H):=\lim_{n\rightarrow\infty}C_n(H)$ is the \emph{modified
Chebyshev constant}.

The following result on the relation of $M$ and $D$ was shown in \cite{Bela}, but follows also
from the combination of \cite{ChoquetSem} and \cite{Oh2}.
\begin{proposition}
\label{Mesdelta} For any $n\in\bN$ and $H\subset X$ we have
$D_n(H)\le M_n(H)$, and thus also $D(H)\le M(H)$.
\end{proposition}
 Let us recall the connection of the transfinite diameter $D$ and the energy $w$ from
\cite{ChoquetSem}, \cite{Bela}.

\begin{theorem}
\label{deltaesi} The following assertions hold for $D$ and $w$.
\begin{enumerate}[a)]
\item \label{DHlewH} $D(H) \le w(H)$ for all $H\subset X$.
\item \label{DKeqwK} $D(K) = w(K)$ for all $K\Subset X$ compact sets.
\item \label{DHDKDWbdk} If the kernel $k$ is finite-valued (i.e., $0\le k <+\infty$), then  $D(H)=w(H)$ for all $H\subset X$.
\end{enumerate}
\end{theorem}

\noindent  The following proposition is trivial from the definitions, but should be stated
explicitly here.
\begin{lemma}\label{monotoneset}
Let $X$ be a locally compact Hausdorff space and $k$ be any nonnegative, l.s.c., symmetric
kernel on $X$ and $L\subset X$ fixed. Then the set functions
\begin{align*}
H\mapsto & \mbox{$u(H)$, $v(H)$, $w(H)$, $q(H,L)$, $q^{\sstar}(H,L)$, $\underline{q}(L,H)$,
$\underline{q}^{\sstar}(L,H)$,} \notag \\
&\mbox{$D_n(H)$, $D(H)$, $M_n(L,H)$ {\rm and}~ $M(L,H)\in\RRb$}
\end{align*}
are non-increasing. Also, the set-to-set functions
$$ H\mapsto \mbox{$R_n(H,L)$, $R(H,L)$ {\rm and} $A(H,L) \subset \overline{\RR}_{+}$}
$$
are non-increasing. On the other hand for fixed $H\subset X$ the functions
$$ L \mapsto \mbox{$q(H,L)$,
$q^{\sstar}(H,L)$, $\underline{q}(L,H)$, $\underline{q}^{\sstar}(L,H)$, $M_n(L,H)$ {\rm and}
$M(L,H)\in\overline{\RR}_{+}$}
$$
and also the set-to-set functions
$$ L \mapsto \mbox{$R_n(H,L)$, $R(H,L)$ {\rm and} $A(H,L) \subset \overline{\RR}_{+}$}
$$
are non-decreasing.
\end{lemma}
\begin{proof}Trivial by observation of the definitions. Allowing
infima or suprema over larger sets results in smaller or larger values, respectively.
\end{proof}

\begin{lemma} \label{iesm}
Let $H,L\subset X$ be arbitrary.  Then we have $\underline{q}(H,L)\leq q(L,H)$.
\end{lemma}
\begin{proof} Let $L'\subseteq L$ compact and $\mu\in\fM_1(H)$, $\nu\in\fM_1(L')$, both compactly supported in $H$ and $L'$, respectively.
By Fubini's theorem we can write
\begin{equation}
\begin{split}
\underline{Q}(\mu;L')&=\inf_{x\in L'} U^\mu(x)=\inf_{x\in
L'}\int_{\supp\mu} k(x,y)\dd\mu(y)\leq
\int_{L'}\int_{\supp\mu} k(x,y)\dd\mu(y)\dd\nu(x)=\\
&=\int_{\supp\mu}\int_{L'} k(x,y)\dd\nu(x)\dd\mu(y)\leq \sup_{y\in H} \int_{L'} k(x,y)\dd\nu(x)=\\
&=\sup_{y\in H}U^\nu(y)=Q(\nu;H)~.
\end{split}
\end{equation}
Taking supremum and infimum over $\mu$ and $\nu$ respectively it follows, also by applying
Lemma \ref{l:q} and \ref{l:-q-}, that
\begin{equation*}
\underline{q}(H,L')\leq q(L',H).
\end{equation*}
So by Lemma \ref{monotoneset} and taking infimum for $L'$ with the use of Lemma \ref{l:q}
again, we obtain
\begin{equation*}
\underline{q}(H,L)\leq \inf_{L'\Subset L}q(L',H)=q(L,H),
\end{equation*}
hence the assertion.
\end{proof}

Lemma \ref{iesm} and the trivial assertions in Lemma \ref{monotoneset} have the following
corollary.
\begin{corollary}\label{qfordul} Let $H\subset
L\subset X$ be arbitrary subsets. Then we have
$\underline{q}(H,L)\leq q(H,L)$.
\end{corollary}
\begin{proof} In view of Lemma \ref{iesm} we have both
$\underline{q}(H,L)\leq q(L,H)$ and $\underline{q}(L,H)\leq q(H,L)$. On the other hand the
monotonicity properties in Lemma \ref{monotoneset}, imply for $H\subset L$ the inequalities
$q(L,H)\leq q(L,L)\leq q(H,L)$ or, alternatively, $\underline{q}(H,L)\leq \underline{q}(L,L)
\leq \underline{q}(L,H)$. On combining the first, resp.~the second set of inequalities, the
assertion follows both ways.
\end{proof}

Now we turn to further relationships between various capacities
and Chebyshev constants. The following lemma is standard, we only
state it for completeness and because we could not find a
standard reference containing exactly the form we need.
\begin{lemma}\label{lem:KMsuru}
The convex combinations of Dirac measures $\delta_x$ concentrated
on points $x$ of a given compact set $K\subset X$ form a
weak$^*$-dense subset of $\fMe(K)$.
\end{lemma}
\begin{proof} Denote $B$ the unit ball in the Banach space $\fM^{\pm}(K)$.
Writing $\one$ for the constant one function on $K$, we have
$\fMe(K)=\{\mu:\:\mu\in B,\: \int_K \one\dd\mu=1\}$, so $\fMe(K)$
is a weak$^*$-closed subset of $B$. Since $B$ is weak$^*$-compact
in view of the Banach--Alaoglu Theorem, $\fMe(K) \subset B$ is
\ws-compact, too.

It is well-known that the set $E$ of extremal points of $\fMe(K)$
consists of the Dirac measures $\delta_x$, $x\in K$ (see, e.g.,
\cite[Proposition 2.1.2, page 52]{nieberg:1991}). Thus by the Kre\u
\i n--Milman theorem (see, e.g., \cite[Sec.~10.4]{schaefer:1980})
the set of convex combinations $\conv E$ of Dirac measures is
weak$^*$-dense in $\fM_1(K)$.
\end{proof}
By the above lemma, given any measure $\mu\in \fMe(K)$ we can approximate it in the
weak$^*$-topology by a finitely supported measure. On the other hand, finitely supported
measures can further be approximated even in a more strict sense.
\begin{lemma}\label{lem:raci}
For any finitely supported probability measure $\nu$ and
$\varepsilon>0$, there exists a probability measure of the form
$\mu=\frac{1}{m}\sum_{i=1}^m \delta_{z_i}$ having the same
support as  $\nu$ and satisfying $(1-\varepsilon)\nu\leq \mu\leq
(1+\varepsilon)\nu$.
\end{lemma}
\begin{proof}
Let $\nu=\frac{1}{n}\sum_{j=1}^n \alpha_j\delta_{w_j}$ with the
coefficients $\alpha_j$ being strictly positive. Choose
$\beta_j\in
((1-\varepsilon)\alpha_j,(1+\varepsilon)\alpha_j)\cap\QQ_+$ such
that $\sum_{j=1}^n \beta_j=1$. Define $\mu=\sum_{j=1}^n
\beta_j\delta_{w_j}$, so for some $m\in \NN$ and with an
appropriate set of points $z_i$ (i.e., with possible repetitions
of the points $w_j$) $\mu$ can be written as
$\mu=\frac{1}{m}\sum_{i=1}^m \delta_{z_i}$.
\end{proof}

The lemma below is a well-known version of the Stone--Weierstra\ss\ Theorem. Actually, without
any assumption on the coefficients  $\alpha_i$, it is the usual approximation theorem used
frequently in the literature. The proof of this slightly ``stronger'' statement goes along the
same lines as that of the original version, see for example Bourbaki \cite[Lemma III.1.2 \&
Exercise III.1.6]{BourbakiXIII}.
\begin{lemma}\label{lem:approxlemma}
Let $f\in C_c(X\times X)$ and $\varepsilon>0$. Then there exists
$g\in C_c(X)\otimes C_c(X)$ of the form
$$
g(x,y)=\sum_{i=1}^n \alpha_i\cdot g_i(x)g'_i(y)
$$
with $g_i,\:g'_i\in C_c(X)$, $0\leq g_i,\:g'_i\leq 1$ and
$|\alpha_i|\leq \|f\|$ satisfying
$$
\|f-g\|\leq \varepsilon.
$$
One can always take also $\alpha_i\in\QQ$. Furthermore, in case
$f\ge 0$ the coefficients $\alpha_i$ can be chosen from $\QQ_+$.
\end{lemma}

\begin{lemma}\label{thm:continuous} Let $K\subset X$ be compact
and $k$ be continuous on ${K\times K}$. Then the mapping
$$
\fM_1(K)\ni\mu\mapsto U^{\mu}\left(=\int_X k(\cdot,y)
\dd\mu(y)\right)
$$
is continuous from the weak$^*$-topology to the sup-norm topology
of $C(K)$.
\end{lemma}
\begin{proof}
Fix $\mu\in \fMe(K)$ and take $\nu\in \fMe(K)$. Further, for
$k|_{K\times K}$ take $g(x,y)=\sum_{i=1}^n \alpha_i\cdot
g_i(x)g'_i(y)$ supplied by Lemma \ref{lem:approxlemma}. Then we
have
\begin{equation*}
\begin{split}
&\int_X k(x,y)\dd \mu(y)-\int_X k(x,y)\dd \nu(y)=\int_X k(x,y)\dd \mu(y)-\int_X g(x,y)\dd \mu(y)+\\
&\qquad+\int_X g(x,y)\dd \mu(y)-\int_X g(x,y)\dd \nu(y)+\int_X
g(x,y)\dd \nu(y)-\int_X k(x,y)\dd \nu(y).
\end{split}
\end{equation*}
So we obtain
\begin{equation}\label{eq:eps3}
\begin{split}
&\Bigl|\int_X k(x,y)\dd \mu(y)-\int_X k(x,y)\dd \nu(y)\Bigr|\leq
\\ &\quad\quad\leq\varepsilon\|\mu\|+\Bigl|\int_X g(x,y)\dd \mu(y)-\int_X
g(x,y)\dd \nu(y)\Bigr|+\varepsilon\|\nu\|.
\end{split}
\end{equation}
Using the particular form of $g$ we can write
\begin{equation*}
\begin{split}
&\Bigl|\int_X g(x,y)\dd \mu(y)-\int_X g(x,y)\dd
\nu(y)\Bigr|=\\
&\qquad\qquad=\Bigl|\sum_{i=1}^n \alpha_i\cdot g_i(x)\Bigl(\int_X
g'_i(y)\dd\mu(y)-\int_X g'_i(y)\dd\nu(y)\Bigr)\Bigr|\leq \\
&\qquad\qquad\leq n\|k\|\max_{i=1,\dots,n} \Bigl|\int_X
g'_i(y)\dd\mu(y)-\int_X g'_i(y)\dd\nu(y)\Bigr|\leq \varepsilon,
\end{split}
\end{equation*}
if $\nu$ is in an appropriate weak$^*$-neighborhood of $\mu$.
Combining this with \eqref{eq:eps3} finishes the proof.
\end{proof}

\begin{theorem}\label{thm:qqstar} Let $k$ be a continuous kernel.
Then $q(H,K)=q^{\sstar}(H,K)$ for all $H\subset X$ and  compact
set $K\subset X$.
\end{theorem}
\begin{proof}
Let us first prove $q(K',K)=q^{\sstar}(K',K)$ whenever $K'\Subset
X$ is compact. From definition $q(K',K) \leq q^{\sstar}(K',K)$ is
trivial. Now if $q(K',K)=+\infty$, there is nothing to prove,
hence let us assume that $q(K',K)<+\infty$. For an arbitrarily
fixed $\varepsilon>0$ take $\mu$ a probability measure supported
in $K'$: $\supp \mu\subset K'$ such that (by \eqref{qcapacity})
$$
Q(\mu;K)\leq q(K',K)+\varepsilon.
$$
Using Lemma \ref{lem:KMsuru} and Lemma \ref{lem:raci} we find
points $x_j\in K'$ ($j=1,\dots,n$, $n\in\NN$) such that the
measure $\nu=\frac{1}{n}\sum_{i=1}^n \delta_{x_i}$ is close to
$\mu$ (lies in an arbitrarily prescribed weak$^*$-neighborhood of
$\mu$), so by Lemma \ref{thm:continuous} the estimate
$|Q(\mu;K)-Q(\nu;K)|\leq\varepsilon$ holds. So $q^{\sstar}(K',K)
\le Q(\nu;K) \leq q(K',K)+2\varepsilon$, and taking infimum in
$\varepsilon$ yields the statement for compact sets $K'$. We
conclude the proof by referring to Lemma \ref{l:q}.
\end{proof}

\section{Dual Chebyshev constants}

Analogously to the dual formulations $q$ and $\underline{q}$, we
define the dual Chebyshev constants. This will enable us to
identify the rendezvous intervals.

\begin{definition}For arbitrary $H,L \subset X$ the
\emph{$n^{\text{th}}$ (general) dual Chebyshev constant} of $L$
relative to $H$ is defined as
\[
    \overline{M}_n(H,L):= \inf_{w_1,\ldots,w_n \in H} \sup_{x\in L}
    {\frac{1}{n}}\Bigg(
   \sum_{j=1}^{n} k(x,w_j)
   \Bigg)
   \ .
\]
\end{definition}
Showing the quasi-monotonicity of $\overline{M}_n(H,L)$ and using
Lemma \ref{l:quasi} give immediately the following.
\begin{proposition}The sequence $\overline{M}_n(H,L)$ of
dual Chebyshev constants has limit.
\end{proposition}
\begin{proof}
Suppose that $\overline{M}_n(H,L)<+\infty$ for some $n\in\mathbb{N}$. This means that there are
$x_1,x_2,\dots,x_n\in H$ with
\begin{equation*}
\sup_{y\in L} \frac{1}{n}\sum_{i=1}^n k(x_i,y)<+\infty
\end{equation*}
But then taking $x_{n+1}:=x_n$, we see that
$\overline{M}_{n+1}(H,L)$ is also finite. Summing up, we see that
the sequence $\overline{M}_n(H,L)$ is either constant $+\infty$ or
eventually finite in which case Fekete's Lemma  \ref{l:quasi} is
applicable, because $\overline{M}_n(H,L)$ is quasi-monotone
decreasing. In both cases the limit exists.
\end{proof}

Note that by quasi-monotonicity we have $\overline{M}_n(H,L)\ge
\overline{M}(H,L)$ for $n\in\NN$. The limit is denoted by
$\overline{M}(H,L)$, and in particular by $\overline{M}(H)$.

It is easy to see that for arbitrary sets $H,L\subset X$ the
$n^{\text{th}}$ Chebyshev constants and dual Chebyshev constants
are just the lower, resp.~upper endpoints of the intervals of the
$n^{\text{th}}$ rendezvous sets \eqref{nthrandi} and rendezvous
sets \eqref{randidef}; and similarly with the energies $q(H,L)$
and dual energies $\underline{q}(H,L)$ regarding \eqref{avridef}.
\begin{proposition}\label{prop:reni} For arbitrary subsets $H, L \subset X$
we have
\begin{align*}
&\inf R_n(H,L)= M_n(H,L),& \qquad\qquad\qquad &\sup
R_n(H,L)=\overline{M}_n(H,L),
\notag \\
&\inf R(H,L)= M(H,L),& &\sup R(H,L)=\overline{M}(H,L),
\\
&\inf A(H,L)=\underline{q}(H,L), & & \sup A(H,L)=q(H,L)~.&&\notag
\end{align*}
That is, we have
\begin{align}\label{randneqchebyn}
&R_n(H,L)= [M_n(H,L),\overline{M}_n(H,L)], \notag \\
&R(H,L)= [M(H,L), \overline{M}(H,L)], \\
&A(H,L)=[\underline{q}(H,L),q(H,L)]~.\notag
\end{align}
\end{proposition}
\begin{proof}
By Remark \ref{AmH} all the intervals are intersections of
certain families $\mathcal{F}$ of closed intervals. Hence the
supremum (resp.~infimum) of the lower (resp.~upper) endpoints of
the elements of $\mathcal{F}$ is the lower (resp.~upper) endpoint
of $\bigcap \mathcal{F}$.
\end{proof}
\begin{remark}\label{rem:convent}
It is important to note that intervals appearing in the
proposition above \emph{may be empty}. This occurs, for example
for $R(H,L)$,  if $\overline{M}(H,L)<M(H,L)$\footnote{Here and
also in the sequel, we use the convention $[a,b]=\emptyset$,
whenever $b<a$.}. Thus, for example, proving that the rendezvous
interval is non-empty is the same as showing
$M(H,L)\leq\overline{M}(H,L)$.
\end{remark}
\begin{proposition}\label{prop:dchc} For any $H,L\subset X$ the dual
Chebyshev constant $\overline{M}(H,L)$ coincides with
$q^{\sstar}(H,L)$. Also, the dual statement
$M(H,L)=\underline{q}^{\sstar}(H,L)$ holds true.
\end{proposition}
\begin{proof}
The inequality $q^{\sstar}(H,L)\leq \overline{M}_n(H,L)$ is
trivial since in the definition of $\overline{M}_n(H,L)$ we allow
only a subclass of finitely supported probability measures. Hence
$q^{\sstar}(H,L)\leq \overline{M}(H,L)$. If $q^{\sstar}=+\infty$,
we are ready. Otherwise fix $\eta>0$ arbitrarily and take
$\nu=\sum_{j=1}^k \alpha_j\delta_{x_j}$, $x_j\in H$ such that
$Q(\nu;L)<q^{\sstar}(H,L)+\eta$. For $\varepsilon>0$ Lemma
\ref{lem:raci} supplies a measure $\mu=\frac{1}{m}\sum_{i=1}^m
\delta_{z_i}\in\fMe(H)$ with
$$
Q(\mu;L)=\sup_{y\in L}\frac{1}{m}\sum_{i=1}^m
k(z_i,y)\leq\sup_{y\in L}\sum_{j=1}^n (1+\varepsilon)\alpha_j
k(x_j,y)\leq (1+\varepsilon)(q^{\sstar}(H,L)+\eta).
$$
Since this holds for all $\varepsilon, \eta >0$, the proof
concludes. The case of $\underline{q}^{\sstar}$ is similar.
\end{proof}

\begin{corollary} For any $H\subset X$ one has $\underline{q}^{\sstar}(H)=M(H)$ and
$q^{\sstar}(H)=\overline{M}(H)$.
\end{corollary}

\section{Chebyshev constant and energy}

The results presented in this section are complementary to those
in Lemma \ref{iesm}. Indeed, the mentioned lemma gives
$\underline{q}(H,L)\leq q(L,H)$, but we prove here equality under
certain assumptions. The relation to $M(H,L)$ is also observed.

Already Fuglede \cite{Fug2} proves (even under more general
assumptions) the relationship between the dual notions of
energies.

\begin{remark}To prove $q(H,L)=\underline{q}(L,H)$ obviously
requires some minimax results. If our kernel was finite-valued, then a perfect tool would be
Kneser's or Kassay and Ko\-lum\-b\'{a}n's minimax theorem. In the setting of rendezvous numbers,
Thomassen applies the Neumann minimax theorem, Stadje refers to a more game-theoretically
formulated minimax theorem in Ferguson \cite{ferguson:1967}.  Why is this confusing abundance
of occurrences of various minimax results? Probably, the reason is best explained by Frenk,
Kassay and Kolumb\'{a}n in \cite{frenk/kassay/kolumban:2004}, who in fact show that all these and
many other well-known minimax results (e.g., Ky-Fan, Kakutani) are (more or less elementarily)
equivalent to each other. However, our kernel may take infinite values, thus we need a
respective minimax theorem. A possible choice would be to refer to Glicksberg's apparently
unfindable work \cite{glicksberg}, which is mentioned by Fuglede \cite{Fug2}. Actually, Fuglede
tackles the problem by an approximation argument. Instead we choose to present our version with
a full proof, which is based on the lecture notes of Pollard \cite{pollard:2001}.
\end{remark}

Let $C$ be a convex set in a linear space. A function $f:C\to (-\infty,+\infty]$ said to be
convex, if
\begin{equation*}
f(\alpha x+(1-\alpha)y)\leq \alpha f(x)+(1-\alpha) f(y),
\end{equation*}
the multiplication $0\cdot f(x)$ is interpreted as $0$ and $\alpha\cdot +\infty=+\infty$,
$\alpha\neq 0$. Analogously one defines concavity.

\begin{theorem}[Minimax Theorem]\label{thm:minimax1}
Let $A$ be a compact, convex subset of a Hausdorff topological
vector space $U$ and $B$ be a convex subset of the linear space $V$.
Let $k: A\times B\to(-\infty,+\infty]$ be l.s.c.~on $A$ for fixed
$y\in B$, and assume that $k$ is convex in the first and concave in
the second variable. Then
\begin{equation*}
\sup_{y\in B} \inf_{x\in A}k(x,y)=\inf_{x\in A}\sup_{y\in B} k(x,y).
\end{equation*}
\end{theorem}
\begin{proof}
The inequality ``$\leq $'' is trivial. For the converse we may assume that the left hand side
is finite, otherwise we are done. Let both $R>S$ be strictly greater than the left hand side.
We have to show that the right hand side does not exceed $R$. For this purpose define
$A_{y,r}:=\{x\in A:\: k(x,y)\leq r\}$ for $y\in B$ and $r\in\RR$. By assumption on $R$ we have
$A_{y,R}\neq \emptyset$ for all $y\in B$. Now we have to prove that also $\bigcap_{y\in
B}A_{y,R}\neq \emptyset$. Indeed, if $x_0$ belongs to this intersection, then $\sup_{y\in B}
k(x_0,y)\leq R$ and the claim follows.

So let us turn to proving the non-emptyness of the intersection. Since the sets $A_{y,R}$ are
all compact (and non-empty), it suffices to prove that every finite family $A_{y_i,R}$,
$i=1,\dots, n$ ($n\in\NN$) has non-empty intersection.

Let first $n=2$. We argue by contradiction, i.e., suppose that $A_{y_1,R}\cap
A_{y_2,R}=\emptyset$. Introduce the notations $A':=\{x:\:x\in A,\:
k(x,y_1)<+\infty,~k(x,y_2)<+\infty\}$ and $A_j:=A_{y_j,R}\cap A'$, $j=1,2$. By assumption, all
these sets are convex, too. We will first find an $\alpha\in[0,1]$ such that
\begin{equation}\label{eq:minimax2}
\alpha k(x,y_1)+(1-\alpha)k(x,y_2)\geq R\qquad(\forall x\in A').
\end{equation}
If $x\not \in A_{y_1,R}\cup A_{y_2,R}$, then the above inequality
holds for all $\alpha\in[0,1]$. In case $x\in A_1$, we must choose
$\alpha$ so that
\begin{equation*}
k(x,y_2)-R \geq \alpha(k(x,y_2)-k(x,y_1))\qquad (\forall x\in
A_1).
\end{equation*}
In view of the disjointness assumption $A_{y_1,R}\cap
A_{y_2,R}=\emptyset$ we have for any $x\in A_{y_1,R}$
$k(x,y_2)>R\geq k(x,y_1)$, hence we get the equivalent form
\begin{equation*}
\big( 1 \geq \big)~ \inf_{x\in A_1}\frac{k(x,y_2)-R}{k(x,y_2)-k(x,y_1)}\geq \alpha.
\end{equation*}
Similarly, for the case $x\in A_2$ the constraint on $\alpha$ will
be $\alpha(k(x,y_1)-k(x,y_2)) \geq R-k(x,y_2)$ ($\forall x\in
A_2$), that is
\begin{equation*}
\alpha\geq \sup_{x\in
A_2}\frac{R-k(x,y_2)}{k(x,y_1)-k(x,y_2)}~\big(\geq 0\big).
\end{equation*}
Note that by the assumption of disjointness, the two constraints are independent, thus $\alpha$
can be chosen satisfying both requirements if and only if
\begin{equation*}
\frac{R-k(x_2,y_2)}{k(x_2,y_1)-k(x_2,y_2)}\leq
\frac{k(x_1,y_2)-R}{k(x_1,y_2)-k(x_1,y_1)} \qquad \big(\forall
x_1\in A_1~~{\rm and}~~ \forall x_2\in A_2\big).
\end{equation*}
As the denominators are positive, some computation leads to the equivalent formulation
\begin{equation}\label{eq:minimax1}
(R-k(x_2,y_2))(R-k(x_1,y_1))\leq (k(x_1,y_2)-R)(k(x_2,y_1)-R)\quad
\big(x_j\in A_j,~ j=1,2\big).
\end{equation}
So let $x_1\in A_1$ and $x_2\in A_2$ be arbitrary. Note that the terms on the right hand side
are necessarily positive in view of the choice of $x_1$, $x_2$ and the assumption of
disjointness. So we may assume that also on the left hand side of \eqref{eq:minimax1} both
terms are strictly positive. Then there exists $\theta\in (0,1)$ with
\begin{equation}\label{eq:theta}
\theta k(x_1,y_1)+(1-\theta)k(x_2,y_1)=R.
\end{equation}
Define $x_\theta=\theta x_1+(1-\theta) x_2$ for this particular $\theta$. By convexity
$k(x_\theta,y_1)\leq R$, so $x_\theta\in A_{y_1,R}$ and therefore $x_\theta\notin A_{y_2,R}$.
On the other hand
\begin{equation}\label{eq:xtheta}
 +\infty>\theta k(x_1,y_2)+(1-\theta)k(x_2,y_2)\geq k(x_\theta,y_2)>R.
\end{equation}
Combining \eqref{eq:theta} and \eqref{eq:xtheta} we arrive at
\begin{align*}
\frac{R-k(x_1,y_1)}{k(x_2,y_1)-R}=\frac{1-\theta}{\theta}<\frac{k(x_1,y_2)-R}{R-k(x_2,y_2)},
\end{align*}
which gives \eqref{eq:minimax1}.

Suppose that incidentally $\alpha=1$. Then
\begin{equation*}
k(x,y_1)\geq R
\end{equation*}
holds for all $x\in A'$. Note that $M:=\inf_{x\in A} k(x,y_2)$ is a
finite minimum value because of the lower semicontinuity of $k$ and
the compactness of $A$. Now let $\beta:=(S-M)/(R-M)$ and $y_0:=\beta
y_1+(1-\beta) y_2$. By concavity of $k$ in the second variable, we
find
\begin{equation*}
k(x,y_0)\geq \beta k(x,y_1)+(1-\beta)k(x,y_2) \geq
\begin{cases}
+\infty, &{\rm if}~x\not\in A', \\
\beta R +(1-\beta)M=S, &{\rm if}~x\in A'.
\end{cases}
\end{equation*}
Hence in both cases $k(x,y_0)\geq S$. It follows that $\inf_{x\in A} k(x,y_0) \geq S$, a
contradiction with the choice of $S$. Therefore, we conclude that $\alpha=1$ is not possible.

Similarly we find that $\alpha= 0$ can not hold. So we obtain that there is $\alpha\in(0,1)$
satisfying \eqref{eq:minimax2}. But then, actually, for all $x\in A$, since $\alpha\in (0,1)$.
By the concavity of $k$ this implies
\begin{equation*}
 \inf_{x\in A}k(x,\alpha y_1+(1-\alpha)y_2)\geq R,
\end{equation*}
 which, in view of $\alpha y_1+(1-\alpha)y_2\in B$ would be a contradiction to the
choice of $R$.

For completing the proof, we apply induction. Let $n\geq 2$ and assume that the assertion has
already been proved for $n$. Then replacing $A$ by the compact, convex set $A_{y_{n+1},R}$ and
considering the sets $A_{y_{m},R}\cap A_{y_{n+1},R}$, $m=1,\dots, n$, we can deduce
$\bigcap_{m=1}^{n+1} A_{y_m,R}\ne \emptyset$, too. This concludes the proof.
\end{proof}
Now we are in the position to investigate the relationship between $q$ and $\underline{q}$.
\begin{theorem}[{\bf Fuglede}]\label{thm:minimax}
Let $k$ be any symmetric, non-negative kernel, $K\Subset X$
compact and $L\subset X$ be any subset. Then the equality
\begin{align*}
q(K,L)&= \inf_{\mu\in\fMe(K)}\sup_{\nu\in\fMe(L)}\int_L\int_K
k(x,y)\dd\mu(y)\dd\nu(x)=\\
&=\sup_{\nu\in\fMe(L)} \inf_{\mu\in\fMe(K)}\int_L\int_K
k(x,y)\dd\mu(y)\dd\nu(x)
\end{align*}
holds true. Furthermore, one has $q(K,L)=\underline{q}(L,K)$.
\end{theorem}
\begin{proof} Obviously
\begin{align}\label{qKH:infsup}
q(K,L)&=\inf_{\mu\in\fMe(K)}\sup_{\nu=\delta_z\atop z\in
L}\int_L\int_K k(x,y)\dd\mu(y)\dd\nu(x)\leq\\\notag &\leq
\inf_{\mu\in\fMe(K)}\sup_{\nu\in\fMe(L)}\int_L\int_K
k(x,y)\dd\mu(y)\dd\nu(x)
\end{align}
holds. Observe that $A:=\fMe(K)$ is a convex, nonempty, and
weak$^*$-compact subset of $U:=\fM^{\pm}(K)$, while $B:=\fMe(L)$
is a convex, nonempty subset of the vector space
$V:=\fM^{\pm}(X)$; moreover, the mapping
$$
f:A\times B \to \RR \qquad f(\mu,\nu):=\int_L\int_K
k(x,y)\dd\mu(y)\dd\nu(x)
$$
is affine on $A\times B$ (in fact linear on $U\times V$), while for any fixed measure $\nu\in
B:=\fMe(L)$ it is lower semicontinuous in $\mu\in A:=\fMe(K)$ by Lemma \ref{lem:lsc} b). Hence
we can continue \eqref{qKH:infsup} by an application of the Minimax Theorem \ref{thm:minimax1}
to obtain
\begin{align}\label{qKH:supinf}
q(K,L)& \le \inf_{\mu\in\fMe(K)} \sup_{\nu\in\fMe(L)}
\int_L\int_K k(x,y)\dd\mu(y)\dd\nu(x)=\notag\\
&= \sup_{\nu\in\fMe(L)} \inf_{\mu\in\fMe(K)}\int_L\int_K
k(x,y)\dd\mu(y)\dd\nu(x)=\\\notag &=\sup_{\nu\in\fMe(L)}
\inf_{\mu\in\fMe(K)}\int_L U^{\mu}(x)\dd\nu(x) \le
\sup_{\nu\in\fMe(L)} \inf_{\mu\in\fMe(K)} \sup_{x\in L}
U^{\mu}(x)=\\ \notag & = \inf_{\mu\in\fMe(K)} Q(\mu;L) = q(K,L) ~.
\end{align}
Thus all expressions in \eqref{qKH:supinf} must be equal.

By Lemma \ref{iesm} and by what we proved in \eqref{qKH:supinf}
\begin{align*}
\underline{q}(L,K)&\leq q(K,L)= \sup_{\nu\in\fMe(L)} \inf_{\mu\in\fMe(K)}\int_L\int_K k(x,y)\dd\mu(y)\dd\nu(x)=\\
&=\sup_{\nu\in\fMe(L)} \inf_{\mu\in\fMe(K)}\int_K
U^\nu(y)\dd\mu(y)\leq\sup_{\nu\in\fMe(L)} \inf_{y\in
K}U^\nu(y)=\\
&=\sup_{\nu\in\fMe(L)} \underline{Q}(\nu;K)= \underline{q}(L,K)~,
\end{align*}
where the last inequality follows because $\delta_y\in \fMe(K)$
whenever $y\in K$. This concludes the proof.
\end{proof}

Among other related things, and even under more general assumptions, Ohtsuka shows also the
following results \cite{Oh2}. Notice that again the symmetry of the kernel is actually not
important.
\begin{theorem}[{\bf Ohtsuka}]\label{thm:Mq}
Let $k$ be any symmetric l.s.c.~kernel, $L$ be an arbitrary and
$K$ be a compact subset of $X$. Then the equality
$M(L,K)=\underline{q}(L,K)$ holds.
\end{theorem}
\begin{proof}
Obviously $M(L,K)\leq\underline{q}(L,K)$ holds. First let us
prove the converse inequality in case $L$ is compact, too. Take
$\nu_0\in\fMe(L)$ arbitrarily. Using Lemma \ref{lem:KMsuru} and
Lemma \ref{lem:raci}, we can select a net of measures
$\nu_\alpha:=\frac{1}{n_\alpha}\sum_{i=1}^{n_\alpha}
\delta^\alpha_{w_i}$, $w_i\in L$ weak$^*$-converging to $\nu_0$.

Then by Lemma \ref{lem:lsc} b), we have for each $\alpha$ a
measure $\mu_\alpha\in\fMe(K)$ such that
$$
\int_L\int_K k(x,y)\dd\mu_\alpha(x)\dd\nu_\alpha(y) =
\inf_{\mu\in\fMe(K)}\int_L\int_K k(x,y)\dd\mu(x)
\dd\nu_\alpha(y)~.
$$
Again by weak$^*$-compactness we may assume that $\mu_\alpha$
converges to some $\mu_0\in\fMe(K)$, and again by lower
semicontinuity we find that
$$
\int_L\int_K k(x,y)\dd\mu_0(x)\dd\nu_0(y)\leq \liminf_{\alpha}
\int_L\int_K k(x,y)\dd\mu_\alpha(x)\dd\nu_\alpha(y)~.
$$
These last two relations imply
\begin{equation*}
\begin{split}
&\inf_{\mu\in\fMe(K)}\int_L\int_K k(x,y)\dd\mu(x)\dd\nu_0(y) \leq
\int_L\int_K
k(x,y)\dd\mu_0(x)\dd\nu_0(y)\leq\\
&\quad\leq \liminf_{\alpha} \inf_{\mu\in\fMe(K)}\int_L\int_K
k(x,y)\dd\mu(x)\dd\nu_\alpha(y)=
\\
&\quad=\liminf_{\alpha}\inf_{\mu\in\fMe(K)}\int_K
\frac{1}{n_\alpha}\sum_{i=1}^{n_\alpha}
k(x,w^\alpha_i)\dd\mu(x)\leq\\
&\quad\leq \liminf_{\alpha} \inf_{\mu=\delta_z\atop z\in K}\int_K
\frac{1}{n_\alpha}\sum_{i=1}^{n_\alpha} k(x,w^\alpha_i)\dd\mu(x)=\\
&\quad=\liminf_{\alpha}\inf_{ z\in
K}\frac{1}{n_\alpha}\sum_{i=1}^ {n_\alpha} k(z,w^\alpha_i)\leq
\liminf_{\alpha}M_{n_\alpha}(L,K) = M(L,K)~.
\end{split}
\end{equation*}
Taking supremum in $\nu_0\in\fMe(L)$ we obtain
\begin{equation}\label{eq:supinf}
\sup_{\nu\in\fMe(L)}\inf_{\mu\in\fMe(K)}\int_L\int_K
k(x,y)\dd\nu(y)\dd\mu(x)\leq M(L,K)~.
\end{equation}
The left-hand side of \eqref{eq:supinf} equals
$q(K,L)=\underline{q}(L,K)$ by Theorem \ref{thm:minimax}, hence
the assertion follows in case $L$ is compact.

Now let $\emptyset\neq L\subset X$ be arbitrary. Let $\nu_0\in
\fMe(L)$ be with $\supp\nu_0\Subset L$. Then we have
$$\underline{Q}(\nu_0;K)\leq \underline{q}(\supp\nu_0,K)= M(\supp\nu_0,K)\leq M(L,K),$$
where we used what we had proved above for the compact case together with Lemma
\ref{monotoneset}. Now applying \eqref{-q-lemma} from Lemma \ref{l:-q-} concludes the proof.
\end{proof}

\section{Rendezvous numbers}

Before drawing some consequences of the above results, let us
summarize them as follows.
\begin{corollary}\label{cor:sum} Let $H, L\subset X$, then
\begin{alignat}{2}\label{eq:sum1}
 M(H,L)=\underline{q}^{\sstar}(H,L)\leq \underline{q}(H,L)\leq &&\:q(L,H)&\leq q^{\sstar}(L,H)=\overline{M}(L,H)~.\\
 \intertext{If $L\subset X$ is compact, then}
M(H,L)=\underline{q}^{\sstar}(H,L)=\underline{q}(H,L)=&&\:q(L,H)&~. \label{eq:sum2} \\
\intertext{If $K\subset X$ is compact and $k$ is continuous, then}
 &&\:q(L,K)&= q^{\sstar}(L,K)=\overline{M}(L,K)~. \label{eq:sum3}
 \end{alignat}
\end{corollary}
Let us shortly comment on the above proved (in)equalities. This is not closely related to
rendezvous numbers, but complements the abstract potential theoretic point of view.
\begin{remark}
We say that the kernel $k$ satisfies the {\em maximum principle} if for every measure $\mu\in
\fMe$
\begin{equation}\label{maxprin}
  U^\mu(x) \le \sup_{y\in \supp \mu} U^\mu(y)
  \qquad \mbox{for all $x \in X$} \ .
\end{equation}
Note that if the kernel $k$ has the maximum principle then $w(H)=v(H)=u(H)$ as discussed in
\cite{Fug}, but since $w(H)\leq q(H)\leq u(H)$, the energy $q(H)$ equals the other ones, too.

In \cite{Bela} it is shown that assuming the maximum principle for the kernel implies
$D(K)=w(K)=u(K)=M(K)$ for all $K\subset X$ compact set. In this case we also have
$w(K)=q(K)=u(K)=M(K)$. Without assuming the maximum principle, in general we only know
$M(K)=q(K)$ by Equality \eqref{eq:sum2} in the above corollary. Keeping in mind Theorem
\ref{deltaesi}, we see that indeed the maximum principle implies the equality $D(K)=M(K)$.
Recall that, in the classical case of the logarithmic kernel $k(x,y)=-\log|x-y|$ on $\CC$, the
equality of the two Chebyshev constants $C(K)=M(K)$ holds. Indeed, as proved above
$M(K)=M(K,K)=q(K)$ and $C(K)=M(\CC,K)=q(K,\CC)=u(K)$. So again the maximum principle is the
reason for $M(K)=C(K)$.
\end{remark}

Let us start with a result showing that our definitions for
rendezvous intervals are non-trivial.
\begin{theorem}\label{th:existence} Let $X$ be a locally compact Hausdorff
space, $\emptyset \ne H \subset L \subset X$ be arbitrary, and let
$k$ be any nonnegative, symmetric kernel on $X$. Then the
intervals $R_n(H,L)$, $R(H,L)$ and $A(H,L)$ are nonempty.
\end{theorem}
\begin{proof}
It suffices to prove the assertion concerning $A(H,L)$ since it
is contained in the other two sets (cf.~Remark \ref{AmH}). By
Proposition \ref{prop:reni}, $A(H,L)=[\underline{q}(H,L),q(H,L)]$
holds, and the assertion reduces to stating
$\underline{q}(H,L)\leq q(H,L)$, which is furnished by Corollary
\ref{qfordul} or \eqref{eq:sum1} (cf.~also Remark
\ref{rem:convent}).
\end{proof}

\begin{remark} It is easy to see that in case $k:\Delta\to\{+\infty\}$
(where $\Delta$ is the diagonal $\{(x,x)~:~x\in X\}$) we have
$q^{\sstar}(H)=+\infty$, hence the rendezvous intervals $R(H)$
always extend to $+\infty$. Note also that in general $q(H)\not=
+\infty$, providing examples where $A(H)\subsetneqq R(H)$, and in
particular, that $R(H)$ is not unique.
\end{remark}

\begin{example}
Take, e.g., the case of $X:=H:=[0,1]$ and $k(x,y):=-\log |x-y|$, where obviously
$\overline{M}(H)=q^{\sstar}(H)=+\infty$, while $M(H)=\log 4$, since we know that classical
capacity of $[0,1]$ is $1/4$ (see, e.g., \cite[Cor.~5.2.4]{Rans}). Further $q(H)=\log 4$ holds
also, because $w([0,1])\leq M([0,1])\leq q([0,1])\leq u([0,1])$, but, as well-known, $k$ has
maximum principle, so $w=q=u$.
\end{example}

\begin{remark} It is easy to construct examples, when the
rendezvous intervals are ``almost empty'': consider, e.g.,
$R_n(\RR,\RR)=\{+\infty\}$. This and Remarks \ref{sdc} and
\ref{sd} explain the slightly disturbing situation that some
papers state that ``there is no rendezvous number'' for cases
where we find one. However, not only $+\infty$ can show up in the
closure of intervals for the definition of rendezvous numbers,
hence not only $+\infty$ can be a rendezvous number for us while
does not exist for other authors. For the case of the $\ell_p$
spaces see \cite{FRnorm}.
\end{remark}

\begin{theorem}\label{th:unique} Let $X$ be any locally compact Hausdorff
topological space, $k$ be any l.s.c., nonnegative, symmetric
kernel function, and $\emptyset\neq K\Subset X$ compact. Then
$A(K)$ consists of one single point. Furthermore, if $k$ is
continuous, then even $R(K)$ consists of only one point.
\end{theorem}
\begin{proof} The first part follows directly from Proposition
\ref{prop:reni}, \eqref{randneqchebyn} and taking $L=K$ in the second part of Theorem
\ref{thm:minimax} or in \eqref{eq:sum2} in Corollary \ref{cor:sum}. For the second assertion,
after taking $L=K$ both \eqref{eq:sum2} and \eqref{eq:sum3} in Corollary \ref{cor:sum} can be
applied, hence the assertion follows.
\end{proof}
\begin{example}
Uniqueness for $R(K,L)$ (or $A(K,L)$) can not be obtained in
general, even if $K\subset L$ and both are compact.  Take, for
example, $X:=\{a,b\}$ a discrete topological space with any
kernel $k$ and $K:=\{a\}$. Then
$$
R(K,X)=[\min \{k(a,a),\: k(a,b)\},\: \max \{k(a,a),\: k(a,b)\}]~,
$$
which reduces to one point only if $k(a,\cdot)$ is constant.
\end{example}

\begin{theorem}\label{thm:AeqR}
If the kernel $k$ is continuous and $L$ is compact we have the
equality $R(H,L)=A(H,L)$ for all $H\subset X$.
\end{theorem}
\begin{proof}
Proposition \ref{prop:reni} identifies the rendezvous intervals
and the results summarized in Corollary \ref{cor:sum} conclude the
proof.
\end{proof}

\begin{remark} The rendezvous set $R(H,L)$ is independent of any
topology, while the notion of $A(H,L)$ essentially refers to
regular Borel measures. Thus in some cases equality of them is
useful.
\end{remark}

\section{Concluding remarks and hints of further work}\label{sec:conclusion}

The modified Definitions \ref{randidef} and \ref{avridef} of
rendezvous numbers and average numbers lead to quite general
existence results, well over the restrictions usual in the
theory. The reason is the use of closure: in many cases, e.g., in
cases when the kernel is not continuous, but only l.s.c.,
considering closure saves the day. Basic results of the theory of
rendezvous numbers extend quite well in the new setting. For
concrete results and also as regards heuristic ideas and general
perception of the topic, general potential theory with a
l.s.c.~kernel on a locally compact, Hausdorff space turns out to
be the relevant setting.

Given the above state of the matter, we aimed at understanding
the working force and the general principles behind further,
particular results. In metric spaces, there is a theory of
rendezvous numbers related to invariant measures \cite{MN} and in
relation to maximal energy \cite{Bj, W3}. These results -- even
if not the available proofs! -- can all be conveniently described
by potential theory, hence it is natural to expect general
versions of the results known so far. For these see \cite{FRmetr}.

As for Banach spaces, extension of the existence and uniqueness results from locally compact
spaces even to infinite dimensional normed spaces, deserves attention. We can accomplish this,
showing that the new definition works definitely better than the old one of, say, strong
rendezvous numbers without use of closure. One can describe a few further, fairly general
results on rendezvous numbers of normed spaces, and therefore computations of rendezvous sets
and numbers of concrete normed spaces become accessible. All this is most interesting in cases
where up to now general understanding stopped at the fact that (\emph{strong}) rendezvous
numbers do not exist. Our results regarding the above questions will be presented in the
forthcoming work \cite{FRnorm}.

\section{Acknowledgements}
For fruitful discussions, comments and references at the outset of this work, we are indebted
to Norm Levenberg. We owe particular gratitude to Natalia Zorii for calling our attention to
\cite{Oh2}, which finally led us to realize that much of what we constructed on our way to
describe the rendezvous numbers, had already been made ready in general linear potential theory
by Fuglede and Ohtsuka. We also thank the anonymous referee for a careful reading and several
valuable suggestions. We strongly hope that the thorough revision due to these references and
suggestions improved the presentation of the material.

\parindent0pt


\begin{thebibliography}{99}

\bibitem{VA}
\textsc{V.~Anagnostopoulos, Sz.~Gy.~R\'{e}v\'{e}sz},
\newblock Polarization constants for products of linear functionals over $\mathbb{R}^2$ and $\mathbb{C}^2$ and Chebyshev constants of the unit sphere,
\newblock \textit{Publ. Math. Debrecen}, \textbf{68} (1--2) (2006), 75--83, {to appear}.

\bibitem{baronti/casini/papini:2000}
\textsc{M.~Baronti, E.~Casini, P.~L.~Papini},
\newblock On average distances and the geometry of Banach spaces,
\newblock \textit{Nonlinear Anal., Theory Methods Appl.} {\bf 42A} (2000), no.~3, 533--541.

\bibitem{BourbakiXIII}
\textsc{N.~Bourbaki},
\newblock \textit{\'{E}l\'{e}ments de Math\'{e}matique, XIII.~Integration} Ch.1,2,3 et 4.,
\newblock Hermann, Paris, 1965.

\bibitem{Bj}
\textsc{G.~Bj\"{o}rck},
\newblock Distributions of positive mass, which maximize a certain generalized energy integral,
\newblock \textit{Ark. Mat.} \textbf{3} (1958), 255--269.

\bibitem{C}
\textsc{L.~Carleson},
\newblock \textit{Selected Problems on Exceptional Sets},
\newblock Van Nostrand Mathematical Studies, No. 13
D. Van Nostrand Co., Inc., Princeton, N.J.-Toronto, Ont.-London 1967.

\bibitem{ChoquetSem}
\textsc{G.~Choquet},
\newblock \textit{Diam\`{e}tre transfini et comparaison de diverses capacit\'{e}s},
\newblock S\'{e}minaire de Th\'{e}orie du Potentiel, Facult\'{e} des Sciences de Paris, 1958/59, 7 pages.

\bibitem{CMY}
\textsc{J.~M.~Cleary, S.~A.~Morris, D.~Yost},
\newblock Numerical geometry -- numbers for shapes,
\newblock \textit{Amer. Math. Monthly} \textbf{93} (1986), 260--275.

\bibitem{Bela}
\textsc{B.~Farkas, B.~Nagy},
\newblock Transfinite diameter, Chebyshev constant, and capacity on locally compact spaces,
\newblock \textit{Alfr\'{e}d R\'{e}nyi Institute preprint series, Hung. Acad. Sci.}, \textbf{7/2004}, 10 pages.

\bibitem{FRmetr}
\textsc{B.~Farkas, Sz.~Gy.~R\'{e}v\'{e}sz},
\newblock Rendezvous numbers of metric spaces -- a potential theoretic approach,
\newblock \textit{Arch. Math.}, to appear.

\bibitem{FRnorm}
\textsc{B.~Farkas, Sz.~Gy.~R\'{e}v\'{e}sz},
\newblock Rendezvous numbers in normed spaces,
\newblock \textit{Bull. Austr. Math. Soc.} \textbf{72} (2005), 423--440, to appear.

\bibitem{ferguson:1967}
\textsc{T.S.~Ferguson},
\newblock \textit{Mathematical Statistics. A Decision Theoretic Approach},
\newblock Probability and Mathematical Statistics, Vol. {\bf 1},
\newblock Academic Press, New York, London, 1967.

\bibitem{Fek}
\textsc{M.~Fekete},
\newblock \"{U}ber die Verteilung der Wurzeln bei gewissen algebraischen Gleichungen mit ganzahligen Koeffizienten,
\newblock \textit{Math. Z.} {\bf17} (1923), 228--249.

\bibitem{frenk/kassay/kolumban:2004}
\textsc{J.~B.~G.~Frenk, G.~Kassay, J.~Kolumb\'{a}n},
\newblock {On equivalent results in minimax theory},
\newblock \textit{European J. Oper. Res.}, \textbf{157}, (2004), no.~1, 46--58.

\bibitem{Fug}
\textsc{B.~Fuglede},
\newblock On the theory of potentials in locally compact spaces,
\newblock \textit{ Acta Math.} \textbf{103} (1960), 139--215.

\bibitem{Fug2}
\textsc{B.~Fuglede},
\newblock Le th\'{e}or\`{e}me du minimax et la th\'{e}orie fine du potentiel,
\newblock \textit{Ann Inst. Fourier} \textbf{15} (1965), 65--87.

\bibitem{GVV}
\textsc{J.~C.~Garc\'{\i}a-V\'{a}zquez, R.~Villa},
\newblock The average distance property of the spaces $\ell_{\infty}^n(\mathbb{C})$ and $\ell_{1}^n(\mathbb{C})$,
\newblock \textit{Arch. Math.} \textbf{76} (2001), 222--230.

\bibitem{Gross}
\textsc{O.~Gross},
\newblock The rendezvous value of a metric space,
\newblock \textit{in: Advances in Game Theory}, \textit{Ann. of Math. Studies}, \textbf{52}, Princeton, 1964, 49--53.

\bibitem{glicksberg}
\textsc{I.~L.~Glicksberg},
\newblock Minimax theorem with upper and lower semi-continuous payoffs,
\newblock \textit{Rand. Corp. Res. Memorandum} RM \textbf{478} (1950).

\bibitem{kneser:1952}
\textsc{H.~Kneser},
\newblock {Sur un th\'{e}or\'{e}me fondamental de la th\'{e}orie des jeux},
\newblock \textit{Comptes Rendus Acad. Sci. Paris} \textbf{234} (1952), 2418--2420.

\bibitem{L}
\textsc{N.~S.~Landkof},
\newblock \textit{Foundations of modern potential theory},
\newblock Die Grundlehren der mathematischen Wissenschaften, Band 180. Springer-Verlag, New York-Heidelberg, 1972.

\bibitem{Lin}
\textsc{P.~K.~Lin},
\newblock The average distance property of Banach spaces,
\newblock \textit{Arch. Math.} \textbf{68} (1997), 496--502.

\bibitem{MN}
\textsc{S.~A.~Morris, P.~Nickolas},
\newblock On the average distance property of compact connected metric spaces,
\newblock \textit{Arch. Math.} \textbf{40}(1983), 459--463.

\bibitem{nieberg:1991}
\textsc{P.~Meyer--Nieberg},
\newblock \textit{Banach lattices},
\newblock Springer--Verlag, 1991.

\bibitem{NY}
\textsc{P.~Nickolas, D.~Yost},
\newblock The average distance property for subsets of euclidean spaces,
\newblock \textit{Arch. Math.} \textbf{50} (1988), 380--384.

\bibitem{Oh}
\textsc{M.~Ohtsuka},
\newblock On potentials in locally compact spaces,
\newblock \textit{J. Sci. Hiroshima Univ. ser A} \textbf{25}, 135--352, 1961.

\bibitem{Oh2}
\textsc{M.~Ohtsuka},
\newblock On various definitions of capacity and related notions,
\newblock \textit{Nagoya Math. J.} \textbf{30} (1967) 121--127.

\bibitem{Oh3}
\textsc{M.~Ohtsuka},
\newblock An application of the minimax theorem to the theory of capacity,
\newblock \textit{J. Sci. Hiroshima Univ. ser A} \textbf{29}, 217--221, 1965.

\bibitem{Oh4}
\textsc{M.~Ohtsuka},
\newblock Generalized capacity and duality theorem in linear programming,
\newblock \textit{J. Sci. Hiroshima Univ. ser A} \textbf{30}, 45--56, 1966.

\bibitem{pollard:2001}
\textsc{D.~Pollard},
\newblock The Minimax Theorem,
\newblock \textit{Unpublished lecture notes}, Paris, 2001. http://www.stat.yale.edu/\~{}pollard/Paris2001/lectures.html

\bibitem{Sz}
\textsc{Gy.~P\'{o}lya, G.~Szeg\H o},
\newblock \"{U}ber den transfiniten Durchmesser (Kapazit\"{a}tskonstante) von ebenen und r\"{a}umlichen Punktmengen,
\newblock \textit{J. Reine Angew. Math.} \textbf{165} (1931), 4-49.

\bibitem{PSz}
\textsc{Gy.~P\'{o}lya, G.~Szeg\H o},
\newblock \textit{Problems and Excercises in Analysis}, vol.~{\bf I},
\newblock Die Grundlehren der matehmatischen Wissenschaften in Einzeldarstellungen, {Bd.~193}, Springer Verlag, 1972.

\bibitem{Rans}
\textsc{T. Ransford},
\newblock \textit{Potential Theory in the Complex Plane},
\newblock London Mathematical Society Student Texts {\bf 28},
\newblock Cambridge University Press, 1994.

\bibitem{RS}
\textsc{Sz.~Gy.~R\'{e}v\'{e}sz, Y.~Sarantopoulos},
\newblock Plank problems, polarization, and Chebyshev constants, \textit{J. Korean
Math. Soc.}, \textbf{41} (2004) no. 1, 157--174.

\bibitem{schaefer:1980}
\textsc{H.~H.~Schaefer},
\newblock \textit{Topological Vector Spaces},
\newblock Graduate Texts In Mathematics, vol.~3, Springer-Verlag, 1980.

\bibitem{Stadje}
\textsc{W.~Stadje},
\newblock A property of compact, connected spaces,
\newblock \textit{Arch. Math.} \textbf{36} (1981), 275--280.

\bibitem{Th}
\textsc{C.~Thomassen},
\newblock The rendezvous number of a symmetric matrix and a compact connected metric space,
\newblock \textit{Amer. Math. Monthly} \textbf{107} (2000), no.~2, 163--166.

\bibitem{W1}
\textsc{R.~Wolf},
\newblock On the average distance property in finite dimensional real Banach spaces,
\newblock \textit{Bull. Austral. Math. Soc.} \textbf{51} (1994), 87--101.

\bibitem{W2}
\textsc{R.~Wolf},
\newblock On the average distance property of spheres in Banach spaces,
\newblock \textit{Arch. Math.} \textbf{62} (1994), 338--344.

\bibitem{W3}
\textsc{R.~Wolf},
\newblock On the average distance property and certain energy integrals,
\newblock \textit{Ark. Mat.} \textbf{35} (1997), 387--400.

\bibitem{YZ}
\textsc{L.~Yang, J.~Z.~Zhang,}
\newblock Average distance constants of some compact convex spaces,
\newblock \textit{L. China Univ. Sci. Tech.} \textbf{17}(1987), no. 1, 17--24.

\bibitem{Z}
\textsc{V.~P.~Zaharjuta},
\newblock Transfinite diameter, Chebishev constants, and capacity for compacta in $\mathbb{C}^n$,
\newblock \textit{Math.\ USSR Sbornik}, vol.~{\bf 25} (1975), no.~3, 350--364 (English translation).

\end{thebibliography}
\end{document}